# A geometric interpretation of Hamilton's Harnack inequality for the Ricci flow [*]

Bennett Chow and Sun-Chin Chu
School of Mathematics
University of Minnesota
Minneapolis, MN 55414[†]

## 1 Introduction

In the paper of Li and Yau [LY], a differential Harnack inequality was proved for the heat equation on a Riemannian manifold. Their technique was based upon the maximum principal which made possible the extension to geometric evolution equations. In particular, Richard Hamilton proved differential Harnack inequalities for the mean curvature flow [H1] and the Ricci flow [H2]. He also extended the result of Li-Yau and proved a matrix Harnack inequality for the heat equation [H3]. Recently, Ben Andrews [A] has proved differential Harnack inequalities for very general curvature flows of hypersurfaces, including anisotropic flows.

The purpose of this paper is to give a geometric interpretation of Hamilton's Harnack inequality for the Ricci flow. We shall show that the Harnack quantity is in fact the curvature of a torsion-free connection compatible with a degenerate metric on space-time. More precisely, let $(M, g(t))$ be a solution to the Ricci flow:

$$\frac{\partial}{\partial t} g_{ij} = -2R_{ij}. \tag{1.1}$$

Define the 3-tensor $P$ by:

$$P_{ijk} = \nabla_i R_{jk} - \nabla_j R_{ik}.$$

Since $P$ is antisymmetric in $i$ and $j$, we may consider $P$ as a section of the bundle $\wedge^2 \otimes \wedge^1$ of 2-forms tensor product 1-forms.

---

Define the symmetric 2-tensor[1] $M$ by:

$$M_{ij} = \Delta R_{ij} - \frac{1}{2}\nabla_i\nabla_j R + 2R_{kijl}R_{kl} - R_{ik}R_{kj}.$$

Here $R_{ijkl} = g_{lm}R^m_{ijk}$, and $M_{ij}$ differs from Hamilton's definition by the omission of the term $\frac{1}{2t}R_{ij}$. Hamilton's Harnack inequality says that if $(M, g(t))$ is a solution to the Ricci flow with semi-positive curvature operator and either $(M, g(t))$ is compact or complete with bounded curvature, then for any 1-form $W_i$ and 2-form $U_{ij}$ we have (Theorem 1.1 of [H3]):

$$\begin{aligned} Z &\stackrel{def}{=} M_{ij}W_iW_j - 2P_{ijk}U_{ij}W_k + R_{ijkl}U_{ij}U_{lk} \\ &\geq -\tfrac{1}{2t}R_{ij}W_iW_j, \end{aligned} \quad (1.2)$$

where our definition differs from Hamilton's by the change of sign of the cross term which may be obtained by replacing $W$ by $-W$. Tracing yields the inequality (Corollary 1.2 of [H3]):

$$\frac{\partial R}{\partial t} + \frac{R}{t} + 2\nabla_i R V_i + 2R_{ij}V_iV_j \geq 0 \quad (1.3)$$

For any 1-form $V_i$. Taking $V = 0$ one obtains

$$\frac{\partial R}{\partial t} + \frac{R}{t} \geq 0.$$

The proof of the Harnack quantity involved applying the heat operator $\frac{\partial}{\partial t} - \Delta$ to $Z$ and suitably specifying the covariant derivatives in space of $W$ and $U$ at a point and extending $W$ and $U$ in time in a way to simplify the computation and allow for an application of the maximum principle. Inequality (1.3) has been exploited by Hamilton [H4] to obtain instantaneous higher derivative estimates and also to prove the "Little Loop Lemma".

In section 14 of [H4], following a remark of Nolan Wallach, Hamilton observed an interesting coincidence. The curvature operator $Rm : \wedge^2 M \to \wedge^2 M$ satisfies the equation[2]

$$\frac{\partial}{\partial t}Rm = \Delta Rm + Rm^2 + Rm^\#, \quad (1.4)$$

where

$$(Rm^\#)_{\alpha\beta} = c_\alpha^{\gamma\,\delta}c_\beta^{\mu\,\nu}Rm_{\gamma\mu}Rm_{\delta\nu}$$

is the square of the curvature operator using the Lie algebra structure constants $c_\alpha^{\gamma\,\delta}$ of $so(n)$ which is isomorphic to the fibers of $\wedge^2 M$ (see also section 2 of [H5]).

---

[1] Here and throughout the paper we use the Einstein summation convention. We also do not always bother to raise indices; repeated lower indices is short hand for contraction with respect to the metric.

[2] This equation holds when the computation is made with respect to a time-dependent orthonormal frame, or equivalently, when the tensor $Rm$ is considered as a system of functions on the orthonormal frame bundle (see section 2 of [H2]).

On the other hand, the Harnack quantity $Z$ may be considered as a symmetric element of
$$(\wedge^2 M \oplus \wedge^1 M) \otimes (\wedge^2 M \oplus \wedge^1 M).$$
Define a Lie bracket on the fibers of $\wedge^2 M \oplus \wedge^1 M$ by
$$[U \oplus W, V \oplus X] = [U, V] \oplus (U \rfloor X - V \rfloor W), \tag{1.5}$$
and a degenerate inner product on $\wedge^2 M \oplus \wedge^1 M$ by
$$\langle U \oplus W, V \oplus X \rangle = \langle U, V \rangle. \tag{1.6}$$
Hamilton observed that if one formally writes the equation (compare with (2.10))
$$\frac{\partial}{\partial t} Z = \Delta Z + Z^2 + Z^{\#} \tag{1.7}$$
using the Lie bracket and inner product defined above, one obtains the equation for $Z$ under the Ricci flow, provided the first derivatives of $U$ and $W$ are prescribed suitably (see formula (2.19)-(2.22)). Here $Z$ is considered as a self-adjoint endomorphism of $\wedge^2 M \oplus \wedge^1 M$. This led Hamilton to write:

> "The geometry would seem to suggest that the Harnack inequality is some sort of jet extension of positive curvature operator on some bundle including translation as well as rotation, and this is somehow all related to solitons where the solution moves by translation."

## 2 Main result

Instead of considering a new bundle on the manifold $M$, we consider the tangent bundle of the space-time manifold $M \times [0, T)$, where $[0, T)$ is the time interval of existence of the solution to the Ricci flow. Given $\tau \in [0, T)$, let $\tilde{M}_\tau = M \times [0, T - \tau)$. Define a degenerate metric $\tilde{g}_\tau$ on cotangent space $T^* \tilde{M}_\tau$ by:
$$\tilde{g}_\tau(x, t) = g^{-1}(x, t + \tau)$$
for $(x, t) \in \tilde{M}_\tau$. Here $g^{-1}$ is the inverse of the metric $g$. In local coordinates $\{x^i\}_{i=1}^n$ on $M$ and $x^0 = t$ on $[0, T)$, we have: $\tilde{g}^{ij} = g^{ij}$, if $1 \leq i, j \leq n$, and $\tilde{g}^{ij} = 0$, if $i = 0$ or $j = 0$. Observe that $\tilde{g}_\tau$ is degenerate in the time direction. This implies that one cannot define the Levi-Civita connection in the usual way. That is, the connections compatible with the metric and torsion free are not unique. However, we shall show that one can define a connection $\tilde{\nabla}_\tau$ compatible with the metric $\tilde{g}_\tau$ such that $(\tilde{g}_\tau, \tilde{\nabla}_\tau)$ is a solution to Ricci flow. We say that the pair $(\tilde{g}_\tau, \tilde{\nabla}_\tau)$ is a solution to the Ricci flow for degenerate metrics if it satisfies the system[3]

$$\begin{aligned}\frac{\partial}{\partial t} \tilde{g}^{ij} &= 2 \tilde{g}^{ik} \tilde{g}^{jl} \tilde{R}_{kl} \\ \frac{\partial}{\partial t} \tilde{\Gamma}^k_{ij} &= -\tilde{g}^{kl} (\tilde{\nabla}_i \tilde{R}_{jl} + \tilde{\nabla}_j \tilde{R}_{il} - \tilde{\nabla}_l \tilde{R}_{ij}),\end{aligned} \tag{2.1}$$

---
[3]The general notion of the Ricci flow for degenerate metrics is due to Hamilton.

where $\tilde{\Gamma}_{ij}^k$ are the Christoffel symbols of the connection $\tilde{\nabla}_\tau$. The curvature tensor is defined in the same way as the Riemann curvature tensor:

$$\tilde{R}_{ijk}^l = \partial_i \tilde{\Gamma}_{jk}^l - \partial_j \tilde{\Gamma}_{ik}^l + \tilde{\Gamma}_{jk}^m \tilde{\Gamma}_{im}^l - \tilde{\Gamma}_{ik}^m \tilde{\Gamma}_{jm}^l$$

and the Ricci tensor by $\tilde{R}_{jk} = \Sigma_{p=0}^n \tilde{R}_{pjk}^p$. The reader may be concerned that the curvature tensor does not satisfy all of the usual symmetries of the Riemann curvature tensor, however this will be true for our choices of the connection.

Define the connection $\tilde{\nabla}_\tau$ on $T\tilde{M}_\tau$ by:

$$\tilde{\Gamma}_{ij}^k(x,t) = \Gamma_{ij}^k(x, t+\tau) \text{ if } 1 \leq i,j,k \leq n, \tag{A1}$$

$$\tilde{\Gamma}_{ij}^0(x,t) = 0 \text{ if } 1 \leq i,j \leq n, \tag{A2}$$

$$\tilde{\Gamma}_{i0}^k(x,t) = \tilde{\Gamma}_{0i}^k(x,t) = -R_i^k(x, t+\tau) \text{ if } 1 \leq i,k \leq n, \tag{A3}$$

$$\tilde{\Gamma}_{00}^k(x,t) = -\tfrac{1}{2}\nabla^k R(x, t+\tau) \text{ if } 1 \leq k \leq n \tag{A4}$$

for $(x,t) \in \tilde{M}_\tau$.

The motivation for the definitions above is as follows. For each $t \in [0, T - \tau)$, $M \times \{t\}$ is a hypersurface in $\tilde{M}_\tau$ with induced metric $g(t+\tau)$. Equation (A1) says that the induced connection from $\tilde{\nabla}_\tau$ on $M \times \{t\}$ is the Levi-Civita connection of $g(t+\tau)$. (A2) says that the second fundamental form is identically zero, i.e., $M \times \{t\}$ is totally geodesic. (A3) then implies that $\tilde{\nabla}_\tau \tilde{g}_\tau = 0$, i.e., the connection is compatible with metric. Finally, (A4) implies that the pair $(\tilde{g}_\tau, \tilde{\nabla}_\tau)$ satisfies the Ricci flow. The way we originally obtained the formula (A4) was by guessing which choice would yield $\tilde{R}m = Z$ and then checking that the Ricci flow was indeed satisfied. We would like to have a better understanding of why this choice works.

**Theorem 2.1** *The pair $(\tilde{g}_\tau, \tilde{\nabla}_\tau)$ is a solution to the Ricci flow for degenerate metrics.*

**Theorem 2.2** *The Riemann curvature tensor $\tilde{R}m_\tau$ of $\tilde{\nabla}_\tau$ is given by:*

$$\tilde{R}_{ijk}^l = R_{ijk}^l \text{ if } 1 \leq i,j,k,l \leq n, \tag{B1}$$

$$\tilde{R}_{ijk}^0 = 0 \text{ if } 1 \leq i,j,k \leq n, \tag{B2}$$

$$\tilde{R}_{ij0}^l = -\nabla_i R_j^l + \nabla_j R_i^l \text{ if } 1 \leq i,j,l \leq n, \tag{B3}$$

$$\tilde{R}_{i0k}^l = -\nabla^l R_{ik} + \nabla_k R_i^l \text{ if } 1 \leq i,k,l \leq n, \tag{B4}$$

$$\tilde{R}_{i00}^l = \Delta R_i^l - \tfrac{1}{2}\nabla_i \nabla^l R + 2g^{lm} R_{pim}^q R_q^p - R_m^l R_i^m \tag{B5}$$

$$= \tfrac{\partial}{\partial t} R_i^l - \tfrac{1}{2}\nabla_i \nabla^l R - R_i^m R_m^l \text{ if } 1 \leq i,l \leq n,$$

*where the quantities on the left hand side are evaluated at $(x,t)$ and the quantities on the right are evaluated at $(x, t+\tau)$.*

Theorem 2.1 and 2.2 are special cases of Theorem 3.5 and 3.1 respectively, whose proofs will be given in section 3.

**Corollary 2.3** *The Riemann curvature tensor $\tilde{R}m$ at $(x,t)$ is the same as the Harnack quantity $Z$ at $(x, t+\tau)$.*

Proof. First we observe that the fiber of the bundle of 2-forms $\wedge^2 \tilde{M}_\tau \cong \wedge^2(T^*M \bigoplus R)$ at $(x,t)$ is isomorphic to the fiber of the bundle $\wedge^2 M \bigoplus \wedge^1 M$ at $x$. This may be seen by taking as a basis $\{\omega^i\}$ for the cotangent space $T^*M$ and noting that $\{\omega^i \wedge \omega^j, \omega^k \wedge dt\}$ is a basis for $\wedge^2 \tilde{M}_\tau$. We shall raise an index and consider the following bundle isomorphic to $\wedge^2 \tilde{M}_\tau$:

$$\wedge^{1,1}\tilde{M}_\tau = \{\alpha = \Sigma_{i,j=1}^n \alpha_i^j dx^i \otimes \frac{\partial}{\partial x^j} + \Sigma_{k=1}^n \alpha_k^0 dx^k \otimes \frac{\partial}{\partial t} | \Sigma_{k=1}^n g^{ik}\alpha_k^j = -\Sigma_{k=1}^n g^{jk}\alpha_k^i\}$$
$$\subset T^*\tilde{M}_\tau \otimes T\tilde{M}_\tau.$$

Raising an index of the Riemann curvature tensor so that it is of type (2,2), we have
$$\tilde{R}m : \wedge^{1,1}\tilde{M}_\tau \to \wedge^{1,1}\tilde{M}_\tau$$
is given by
$$\tilde{R}m(T)_j^i = \Sigma_{k,l=0}^n \tilde{R}_{jk}^{il} T_l^k$$
where $\tilde{R}_{jk}^{il} = \tilde{g}^{ip}\tilde{R}_{pjk}^l$ and $T = \Sigma_{p,q=0}^n T_p^q dx^p \otimes \frac{\partial}{\partial x^q} \in \wedge^{1,1}\tilde{M}_\tau$. By the isomorphism $\wedge^{1,1}\tilde{M}_\tau \cong \wedge^{1,1}M \oplus \wedge^1 M$, we may write $T = U \oplus W$, where $U \in \wedge^{1,1}M$ and $W \in \wedge^1 M$. Using Theorem 2.1, we compute

$$\begin{aligned}\tilde{R}m(T,T) &= \tilde{R}_{jk}^{il} U_i^j U_l^k + \tilde{R}_{j0}^{il} U_i^j W_l + \tilde{R}_{0k}^{il} W_i U_l^k + \tilde{R}_{00}^{il} W_i W_l \\ &= R_{jk}^{il} U_i^j U_l^k + (-\nabla^i R_j^l + \nabla_j R^{il}) U_i^j W_l + (-\nabla^l R_k^i + \nabla_k R^{il}) W_i U_l^k \\ &\quad + (\Delta R^{il} - \tfrac{1}{2}\nabla^i \nabla^l R + 2g^{ir} g^{lm} R_{prm}^q R_q^p - R^{im} R_m^l) W_i W_l.\end{aligned} \quad (2.2)$$

Hence by definition (1.2),
$$\tilde{R}m(U \oplus W, U \oplus W) = Z.$$

**Corollary 2.4** *The Ricci tensor $\tilde{R}ic_\tau$ of $\tilde{\nabla}_\tau$ is given by*

$$\tilde{R}_{ij} = R_{ij} \text{ if } 1 \leq i,j \leq n, \tag{C1}$$
$$\tilde{R}_{0j} = \tfrac{1}{2}\nabla_j R \text{ if } 1 \leq j \leq n, \tag{C2}$$
$$\tilde{R}_{00} = \tfrac{1}{2}\frac{\partial R}{\partial t} = \tfrac{1}{2}\Delta R + |Ric|^2 \tag{C3}$$

*where as before the left are at $(x,t)$ and on the right at $(x, t+\tau)$. The scalar curvature $\tilde{R}_\tau$ of $\tilde{\nabla}_\tau$ is given by $\tilde{R}_\tau(x,t) = R(x, t+\tau)$.*

If $V \in TM$, then
$$\tilde{R}ic(V \oplus \frac{\partial}{\partial t}, V \oplus \frac{\partial}{\partial t}) = \frac{1}{2}\frac{\partial R}{\partial t} + \nabla_j R V^j + R_{ij} V^i V^j.$$

The Harnack inequality (1.3), then says that
$$\tilde{R}ic(V \oplus \frac{\partial}{\partial t}, V \oplus \frac{\partial}{\partial t}) \geq -\frac{R}{t}.$$

Note that one can derive the evolution equation for $R$ under the Ricci flow from the second contracted Bianchi identity on $\tilde{M}_\tau$:

$$\frac{1}{2}\frac{\partial R}{\partial t} = \frac{1}{2}\tilde{\nabla}_0 \tilde{R} = \tilde{g}^{ij}\tilde{\nabla}_i \tilde{R}_{j0} = g^{ij}\nabla_i(\frac{1}{2}\nabla_j R) - g^{ij}\tilde{\Gamma}^m_{i0}\tilde{R}_{jm} = \frac{1}{2}\Delta R + |Ric|^2.$$

Similarly, the second Bianchi identity on $\tilde{M}_\tau$ implies the evolution equation (1.4) for $Rm$ by computing

$$\tilde{\nabla}_0 \tilde{R}^l_{ijk} = -\tilde{\nabla}_i \tilde{R}^l_{j0k} - \tilde{\nabla}_j \tilde{R}^l_{0ik}.$$

Observe that the degenerate metric $\tilde{g}_\tau$ defines a metric and a bracket on $\wedge^2 \tilde{M}_\tau$ by

$$\langle \alpha, \beta \rangle = \tilde{g}^{ik}\tilde{g}^{jl}\alpha_{ij}\beta_{kl} \tag{2.3}$$

and

$$[\alpha, \beta]_{ij} = \alpha_{ik}\tilde{g}^{kl}\beta_{lj} - \beta_{ik}\tilde{g}^{kl}\alpha_{lj}. \tag{2.4}$$

Using the decomposition $\wedge^2 \tilde{M}_\tau \cong \wedge^2 M \oplus \wedge^1 M$, we find that the definitions (2.3) and (2.4) above agree with (1.6) and (1.5). Alternately, one may define a metric and a bracket on $\wedge^{1,1}\tilde{M}_\tau$ by

$$\langle \alpha, \beta \rangle = -\alpha^j_i \beta^i_j \tag{2.5}$$

and

$$[\alpha, \beta]^j_i = \alpha^k_i \beta^j_k - \beta^k_i \alpha^j_k. \tag{2.6}$$

This yields equivalent definitions via the isomorphism $\wedge^{1,1}\tilde{M}_\tau \cong \wedge^2 \tilde{M}_\tau$.

Since the pair $(\tilde{g}_\tau, \tilde{\nabla}_\tau)$ satisfies the Ricci flow, the Riemann curvature tensor satisfies the equation[4]:

$$\frac{\partial}{\partial t}\tilde{R}m = \tilde{\Delta}\tilde{R}m + \tilde{R}m^2 + \tilde{R}m^\#, \tag{2.7}$$

where $\tilde{R}m^\#$ is the square of $\tilde{R}m$ using the structure constants of the Lie algebra defined by (2.6). This explains why equation (1.7) holds.

Finally, we compute the evolution equation for $Z = \tilde{R}m(U \oplus W, U \oplus W)$. Let $T = U \oplus W$, from (1.4) we compute:

$$\frac{\partial}{\partial t}Z = \tilde{\Delta}Z + \tilde{R}m^2(T,T) + \tilde{R}m^\#(T,T) + 2\tilde{R}m(\frac{\partial}{\partial t}T - \tilde{\Delta}T, T) - 2\tilde{g}^{ij}\tilde{R}m(\tilde{\nabla}_i T, \tilde{\nabla}_j T)$$
$$- 4\tilde{g}^{ij}\tilde{\nabla}_i \tilde{R}m(\tilde{\nabla}_j T, T).$$

Hence, at a point where

$$\frac{\partial}{\partial t}T = \tilde{\Delta}T \tag{2.8}$$

---

[4] We leave it to reader to check that the computations of the equation for the Riemann curvature tensor in [H5] and [H6] hold for solutions to the degenerate Ricci flow. Again, the equation (2.7) holds with respect to an orthonormal frame in space together with $\frac{\partial}{\partial x^0}$ (since $\tilde{g}_\tau$ is degenerate, this extended frame is not orthonormal).

and
$$\tilde{\nabla}_i T = 0 \text{ for } 1 \leq i \leq n, \tag{2.9}$$

we have
$$\frac{\partial}{\partial t} Z = \Delta Z + \tilde{R}m^2(T,T) + \tilde{R}m^\#(T,T). \tag{2.10}$$

Note that since $Z$ is a scalar, $\tilde{\Delta} Z = \Delta Z$.

In terms of $U$ and $W$, equation (2.8) and (2.9) may be rewritten as:

$$\frac{\partial}{\partial t} U_{ij} = \Delta U_{ij} \tag{2.11}$$

$$\frac{\partial}{\partial t} W_i = \Delta W_i + \frac{1}{2}\nabla^p R\, U_{ip} \tag{2.12}$$

and
$$\nabla_i U_{jk} = 0 \tag{2.13}$$
$$\nabla_i W_j + R_i^p U_{jp} = 0. \tag{2.14}$$

Given equation (2.11) and (2.13), we can relate equations (2.12) and (2.14) to the Ricci solitons. Recall that the Ricci gradient soliton equation is given by (see section 3 of [H2]):
$$R_{ij} = \nabla_i V_j \tag{2.15}$$

where $\nabla_i V_j = \nabla_j V_i$, i.e., $dV = 0$. When $H^1(M;R) = 0$, we have $V = df$ for some function $f$. Tracing (2.15) we obtain:
$$R = \text{div}\, V$$

whereas taking the divergence of (2.15) we have:
$$\frac{1}{2}\nabla_j R + R_{jk} V_k = 0. \tag{2.16}$$

Taking the time derivative of (2.15) and using the equations
$$\frac{\partial}{\partial t} \Gamma^k_{ij} = -g^{kl}(\nabla_i R_{jl} + \nabla_j R_{il} - \nabla_l R_{ij})$$

and
$$\frac{\partial}{\partial t} R_{ij} = \Delta R_{ij} + 2 R_{kijl} R_{kl} - 2 R^2_{ij}$$

we obtain:
$$\nabla_i (\frac{\partial}{\partial t} V_j + R_{jl} V_l) = \nabla_i (\Delta V_j).$$

Assuming there are no parallel 1-forms on $M$, e.g., if $H^1(M;R) = 0$, we have:
$$\frac{\partial}{\partial t} V_j = \Delta V_j - R_{jl} V_l, \tag{2.17}$$

or equivalently,
$$\frac{\partial}{\partial t}V_j = \Delta_d V_j \stackrel{def}{=} -(d\delta + \delta d)V_j. \tag{2.18}$$

Given a 2-form $U_{ij}$ on $M$, define the 1-form
$$W = (U \rfloor V)_i = -U_{ij}g^{jk}V_k.$$

Assuming (2.11), (2.13), (2.15) and (2.17) we compute:
$$\nabla_k W_i = -\nabla_k U_{ij} V_j - U_{ij} \nabla_k V_j = -U_{ij} R_{kj}$$

which is the same as (2.14), and
$$\frac{\partial}{\partial t}W_i = -\frac{\partial}{\partial t}U_{ij} V_j - U_{ij}\frac{\partial}{\partial t}V_j - U_{ij}\frac{\partial}{\partial t}g^{jk}V_k$$
$$= \Delta W_i - U_{ij} R_{jk} V_k,$$

which is equivalent to (2.12) using (2.16).

The computations above are related to section 3 and Theorem 4.1 of Hamilton as follows. Hamilton considers $W$ as the basic quantity and defines $U = V \wedge W$. From the equations[5]

$$\frac{\partial}{\partial t}W_i = \Delta W_i \tag{2.19}$$

and
$$\nabla_i W_j = 0 \tag{2.20}$$

and the Ricci soliton equations, he derives the equations

$$\frac{\partial}{\partial t}U_{ij} = \Delta U_{ij} \tag{2.21}$$

and
$$\nabla_i U_{jk} = \frac{1}{2}(R_{ij}W_k - R_{ik}W_j). \tag{2.22}$$

On the other hand, we consider $U$ as the basic quantity and define $W = U \rfloor V$. From (2.11), (2.13) and the Ricci soliton equations, we derive (2.12) and (2.14). The interesting coincidence to note is that equations (2.11)-(2.14) are equivalent to (2.8) and (2.9). This explains why assuming (2.11)-(2.14) simplifies the computation for $\partial Z/\partial t$. Since (2.19)-(2.22) are related to solitons in an analogous way as (2.11)-(2.14), this gives an additional explanation why assuming (2.19)-(2.22) simplifies the computation in Theorem 4.1 of Hamilton's paper.

---

[5]This equation differs from Hamilton's by $\frac{W_i}{t}$ since we consider steady solitons instead of expanding solitons as our motivation.

## 3  An extension and Proofs

In this section[6] we extend the results of the previous section to the Ricci flow modified by diffeomorphisms generated by gradient vector fields. In particular, given a function $f : M \times [0, T) \to R$, we consider the modified Ricci flow:

$$\frac{\partial}{\partial t} g_{ij} = -2R_{ij} + 2\nabla_i \nabla_j f . \tag{3.1}$$

Here the Hessian of $f$ is the same as the Lie derivative of the metric with respect to the vector field $\nabla f$ i.e., $2\nabla_i \nabla_j f = (\mathcal{L}_{\nabla f} g)_{ij}$.

Analogous to section 1, we consider the degenerate Riemannian manifold $(\tilde{M}_\tau, \tilde{g}_\tau)$. We then look for a connection $\tilde{\nabla}_\tau$ compatible with $\tilde{g}_\tau$ such that the pair $(\tilde{g}_\tau, \tilde{\nabla}_\tau)$ is a solution to the modified Ricci flow for degenerate metrics:

$$\frac{\partial}{\partial t} \tilde{g}^{ij} = 2\tilde{g}^{ik} \tilde{g}^{jl} (\tilde{R}_{kl} - \tilde{\nabla}_k \tilde{\nabla}_l f) \tag{3.2}$$

$$\frac{\partial}{\partial t} \tilde{\Gamma}^k_{ij} = -\tilde{g}^{kl} [\tilde{\nabla}_i (\tilde{R}_{jl} - \tilde{\nabla}_j \tilde{\nabla}_l f) + \tilde{\nabla}_j (\tilde{R}_{il} - \tilde{\nabla}_i \tilde{\nabla}_l f) - \tilde{\nabla}_l (\tilde{R}_{ij} - \tilde{\nabla}_i \tilde{\nabla}_j f)]. \tag{3.3}$$

Note that the modified Ricci flow differs from the Ricci flow by the Lie derivative of $g$ with respect to $\nabla f$ whereas the modified degenerate Ricci flow differs from the degenerate Ricci flow by the Lie derivative of $\tilde{g}_\tau$ with respect to $\tilde{\nabla} f = (\nabla f, \frac{\partial f}{\partial t})$.

The appropriate choice of the connection $\tilde{\nabla}_\tau$ is given by:

$$\tilde{\Gamma}^k_{ij} = \Gamma^k_{ij} \text{ if } 1 \leq i,j,k \leq n, \tag{D1}$$

$$\tilde{\Gamma}^0_{ij} = 0 \text{ if } 0 \leq i,j \leq n, \tag{D2}$$

$$\tilde{\Gamma}^k_{i0} = \tilde{\Gamma}^k_{0i} = -R^k_i + \nabla_i \nabla^k f \text{ if } 1 \leq i,k \leq n, \tag{D3}$$

$$\tilde{\Gamma}^k_{00} = \nabla^k (-\tfrac{1}{2} R + \tfrac{\partial}{\partial t} f - \tfrac{1}{2} |\nabla f|^2 ). \tag{D4}$$

Observe that when $f \equiv 0$, definitions (D1)-(D4) agree with (A1)-(A4).

---

[6] We expect that the results of this section hold when $\nabla_i f$ is replaced by a closed 1-form $V_i$ whose cohomology class is independent of time and $\tilde{\nabla} f$ is replaced by $\tilde{V} = (V, h)$ where $\frac{\partial}{\partial t} V = dh$. In this case one should define $\tilde{\Gamma}^k_{i0} = \tilde{\Gamma}^k_{0i} = -R^k_i + \nabla_i V^k$ and $\tilde{\Gamma}^k_{00} = \nabla^k (-\tfrac{1}{2} R + h - \tfrac{1}{2} |V|^2)$.

**Theorem 3.1** *The Riemann curvature tensor $\tilde{R}m_\tau$ is given by:*

$$\tilde{R}^l_{ijk} = R^l_{ijk} \text{ if } 1 \leq i,j,k,l \leq n, \tag{E1}$$

$$\tilde{R}^0_{ijk} = 0 \text{ if } 1 \leq i,j,k \leq n, \tag{E2}$$

$$\tilde{R}^l_{ij0} = -\nabla_i R^l_j + \nabla_j R^l_i + R^l_{ijp}\nabla^p f \text{ if } 1 \leq i,j,l \leq n, \tag{E3}$$

$$\tilde{R}^l_{i0k} = -\nabla^l R_{ik} + \nabla_k R^l_i + R^l_{ipk}\nabla^p f \text{ if } 1 \leq i,k,l \leq n, \tag{E3'}$$

$$\begin{aligned}\tilde{R}^l_{i00} &= \tfrac{\partial}{\partial t}R^l_i - \nabla^p f \nabla_p R^l_i - \tfrac{1}{2}\nabla_i \nabla^l R - R^p_i R^l_p - \{\nabla_i R^l_p - \nabla_p R^l_i\}\nabla^p f \\
&\quad -\{\nabla^l R_{ip} - \nabla_p R^l_i\}\nabla^p f + R^p_i \nabla_p \nabla^l f - R^l_p \nabla_i \nabla^p f + R^l_{ipq}\nabla^p f \nabla^q f \\
&= \Delta R^l_i - \tfrac{1}{2}\nabla_i \nabla^l R + 2g^{lm}R^q_{pim}R^p_q - R^m_i R^l_m - g^{ql}P_{ipq}\nabla^p f \\
&\quad -g^{ql}P_{qpi}\nabla^p f + R^l_{ipq}\nabla^p f \nabla^q f \text{ if } 1 \leq i,l \leq n.\end{aligned} \tag{E4}$$

Proof. Formulas (E1) and (E2) follow from the definition of $\tilde{R}m_\tau$ and (D1)-(D2). For part (E3) we compute:

$$\begin{aligned}\tilde{R}^l_{ij0} &= \partial_i \tilde{\Gamma}^l_{j0} - \partial_j \tilde{\Gamma}^l_{i0} + \tilde{\Gamma}^p_{j0}\tilde{\Gamma}^l_{ip} - \tilde{\Gamma}^p_{i0}\tilde{\Gamma}^l_{jp} \\
&= \nabla_i(-R^l_j + \nabla_j \nabla^l f) - \nabla_j(-R^l_i + \nabla_i \nabla^l f)\end{aligned}$$

and (E3) follows from the commutation formula

$$\nabla_i \nabla_j \nabla^l f - \nabla_j \nabla_i \nabla^l f = R^l_{ijp}\nabla^p f.$$

Part (E3') is proved similarly. For part (E4) we compute

$$\begin{aligned}\tilde{R}^l_{i00} &= \partial_i \tilde{\Gamma}^l_{00} - \partial_0 \tilde{\Gamma}^l_{i0} + \tilde{\Gamma}^p_{00}\tilde{\Gamma}^l_{ip} - \tilde{\Gamma}^p_{i0}\tilde{\Gamma}^l_{0p} \\
&= \nabla_i \nabla^l(-\tfrac{1}{2}R + \tfrac{\partial f}{\partial t} - \tfrac{1}{2}|\nabla f|^2) + \tfrac{\partial}{\partial t}(R^l_i - \nabla_i \nabla^l f) \\
&\quad - (R^p_i - \nabla_i \nabla^p f)(R^l_p - \nabla_p \nabla^l f)\end{aligned}$$

and (E4) follows from expanding and cancelling off terms.

**Corollary 3.2** *The Ricci tensor $\tilde{Ric}_\tau$ is given by:*

$$\tilde{R}_{ij} = R_{ij} \text{ if } 1 \leq i,j \leq n, \tag{F1}$$

$$\tilde{R}_{0j} = \tfrac{1}{2}\nabla_j R + R_{jp}\nabla^p f \text{ if } 1 \leq j \leq n, \tag{F2}$$

$$\begin{aligned}\tilde{R}_{00} &= \tfrac{1}{2}\tfrac{\partial R}{\partial t} + R_{pq}\nabla^p f \nabla^q f + \nabla^p R \nabla_p f \\
&= \tfrac{1}{2}\Delta R + |Ric|^2 + \nabla^p R \nabla_p f + R_{pq}\nabla^p f \nabla^q f.\end{aligned} \tag{F3}$$

Proof. Apply Theorem 3.1 to the definition $\tilde{R}_{ij} = \tilde{R}^p_{pij}$ while using the contracted second Bianchi identity $\nabla_p R^p_i = \tfrac{1}{2}\nabla_i R$.

An interesting set of identities are the following.

**Lemma 3.3**  
1) $\tilde{R}^k_{i0j} = \tilde{\nabla}_j \tilde{R}^k_i - \tilde{\nabla}^k \tilde{R}_{ij} + \tilde{R}^k_{ipj}\nabla^p f$,  
2) $\tilde{R}^k_{i00} = \tilde{\nabla}_0 \tilde{R}^k_i - \tilde{\nabla}^k \tilde{R}_{0i} + \tilde{R}^k_{ip0}\nabla^p f$,  
3) $\tilde{R}^k_{00j} = \tilde{\nabla}_j \tilde{R}^k_0 - \tilde{\nabla}^k \tilde{R}_{0j} + \tilde{R}^k_{0pj}\nabla^p f = 0$,  
4) $\tilde{R}^k_{000} = \tilde{\nabla}_0 \tilde{R}^k_0 - \tilde{\nabla}^k \tilde{R}_{00} + \tilde{R}^k_{0p0}\nabla^p f = 0$.

Proof. Using the definitions (D1)-(D4) to compare $\tilde{\nabla}_i$ to $\nabla_i$ and $\tilde{\nabla}_0$ to $\frac{\partial}{\partial t}$, we obtain the following. Part 1) follows directly from (E3′) using (E1) and (F1). Part 2) may be derived from (E4) using (F1),(F2) and (E3). Part 3) follows from (F3) and (E3′). Part 4) follows from (F2), (F3) and (F4).

**Remark 3.4** *Formally, parts 2)-4) follows from 1) by setting the appropriate indices to be zero.*

For some of our later computations, it will be convenient to reformulate Lemma 3.3 as follows.

**Corollary 3.5**

1) $\tilde{\nabla}_j(\tilde{R}_{ki} - \tilde{\nabla}_k \tilde{\nabla}_i f) = \tilde{\nabla}_k(\tilde{R}_{ji} - \tilde{\nabla}_j \tilde{\nabla}_i f) + g_{pi}\tilde{R}^p_{kj0}$,  
2) $\tilde{\nabla}_0(\tilde{R}_{ki} - \tilde{\nabla}_k \tilde{\nabla}_i f) = \tilde{\nabla}_k(\tilde{R}_{0i} - \tilde{\nabla}_0 \tilde{\nabla}_i f) + g_{pk}\tilde{R}^p_{i00}$,  
3) $\tilde{\nabla}_j(\tilde{R}_{0k} - \tilde{\nabla}_0 \tilde{\nabla}_k f) = \tilde{\nabla}_k(\tilde{R}_{0j} - \tilde{\nabla}_0 \tilde{\nabla}_j f)$,  
4) $\tilde{\nabla}_0(\tilde{R}_{0k} - \tilde{\nabla}_0 \tilde{\nabla}_k f) = \tilde{\nabla}_k(\tilde{R}_{00} - \tilde{\nabla}_0 \tilde{\nabla}_0 f)$,

where $1 \leq i, j, k \leq n$.

Proof. The corollary follows from Lemma 3.3 and the formula for commuting derivatives:
$$\tilde{\nabla}_i \tilde{\nabla}_j \tilde{\nabla}_k f - \tilde{\nabla}_j \tilde{\nabla}_i \tilde{\nabla}_k f = -\tilde{R}^p_{ijk}\nabla_p f$$
where $0 \leq i, j, k \leq n$.

**Remark 3.6** *Parts 1)-4) are formally equivalent to*
$$\tilde{\nabla}_j(\tilde{R}_{ki} - \tilde{\nabla}_k \tilde{\nabla}_i f) = \tilde{\nabla}_k(\tilde{R}_{ji} - \tilde{\nabla}_j \tilde{\nabla}_i f) + \tilde{R}_{kj0i}$$
*where $0 \leq i, j, k \leq n$ and*
$$\tilde{R}_{kj0i} = \begin{cases} g_{pi}\tilde{R}^p_{kj0} & \text{if } 1 \leq i \leq n \\ 0 & \text{if } i = 0. \end{cases}$$

**Theorem 3.7** *If $g$ is a solution to the modified Ricci flow, then the pair $(\tilde{g}_\tau, \tilde{\nabla}_\tau)$ is a solution to the modified Ricci flow for degenerate metrics.*

Proof. Since $\tilde{g}^{0j} = 0$ for all $0 \leq j \leq n$, to obtain (3.2), we only need to show that
$$\tilde{R}_{kl} - \tilde{\nabla}_k \tilde{\nabla}_l f = R_{kl} - \nabla_k \nabla_l f \qquad (3.4)$$
for $1 \leq k, l \leq n$. This follows from $\tilde{\Gamma}^k_{ij} = \Gamma^k_{ij}$, $\tilde{R}_{ij} = R_{ij}$ (equations (D1) and (F1), respectively). Equation (F1) also implies (3.3) holds for $1 \leq i, j, k \leq n$. Hence we only need to show (3.3) holds in two cases:

**Case 1** $j = 0$ and $1 \leq i, k \leq n$.

The left side of equation (3.3) is:
$$\frac{\partial}{\partial t}\tilde{\Gamma}^k_{i0} = -\frac{\partial}{\partial t}(R^k_i - \nabla_i \nabla^k f).$$

The right side is given by:
$$-\tilde{g}^{kl}\{\tilde{\nabla}_0(\tilde{R}_{il} - \tilde{\nabla}_i\tilde{\nabla}_l f) + \tilde{\nabla}_i(\tilde{R}_{0l} - \tilde{\nabla}_l\tilde{\nabla}_0 f) - \tilde{\nabla}_l(\tilde{R}_{0i} - \tilde{\nabla}_i\tilde{\nabla}_0 f)\}$$
$$= -\frac{\partial}{\partial t}(R^k_l - \nabla_i \nabla^k f) + \tilde{\Gamma}^m_{0i}(\tilde{R}^k_m - \nabla_m\nabla^k f) - \tilde{\Gamma}^k_{0m}(R^m_i - \nabla_i\nabla^m f)$$
$$- g^{kl}\{\nabla_i(\frac{1}{2}\nabla_l R + R^p_l \nabla_p f) - \nabla_l(\frac{1}{2}\nabla_i R + R^p_i \nabla_p f) - \tilde{\Gamma}^p_{i0}R_{pl} + \tilde{\Gamma}^p_{l0}R_{pi} + \tilde{R}^p_{il0}\nabla_p f)\}$$
$$= -\frac{\partial}{\partial t}(R^k_i - \nabla_i\nabla^k f) - g^{kl}(\nabla_i R^p_l - \nabla_l R^p_i + \tilde{R}^p_{il0})\nabla_p f.$$

Equation (3.3) now follows from formula (E3) and the fact that $R^p_{ilq}\nabla^q f \nabla_p f = 0$.

**Case 2** $i = j = 0$ and $1 \leq k \leq n$.

The left side of (3.3) is given by:
$$\frac{\partial}{\partial t}\tilde{\Gamma}^k_{00} = \frac{\partial}{\partial t}[\nabla^k(-\frac{1}{2}R + \frac{\partial f}{\partial t} - \frac{1}{2}|\nabla f|^2)].$$

The right side is:
$$\tilde{g}^{kl}\{-2\tilde{\nabla}_0(\tilde{R}_{0l} - \tilde{\nabla}_0\tilde{\nabla}_l f) + \tilde{\nabla}_l(\tilde{R}_{00} - \tilde{\nabla}_0\tilde{\nabla}_0 f)\}$$
$$= -\tilde{g}^{kl}\tilde{\nabla}_0(\tilde{R}_{0l} - \tilde{\nabla}_0\tilde{\nabla}_l f)$$
$$= -\frac{\partial}{\partial t}[\nabla^k(\frac{1}{2}R - \frac{\partial f}{\partial t} + \frac{1}{2}|\nabla f|^2)]$$
$$+ \tilde{\Gamma}^p_{00}(\tilde{R}^k_p - \tilde{\nabla}_p\tilde{\nabla}^k f) - \tilde{\Gamma}^k_{0p}(\tilde{R}^p_0 - \tilde{\nabla}_0\tilde{\nabla}^p f)$$
$$= \frac{\partial}{\partial t}[\nabla^k(-\frac{1}{2}R + \frac{\partial f}{\partial t} - \frac{1}{2}|\nabla f|^2)],$$

where we used Corollary 3.4 to obtain the first equality and used the identities
$$\tilde{\Gamma}^k_{0p} = -(\tilde{R}^k_p - \tilde{\nabla}_p\tilde{\nabla}^k f) = -(R^k_p - \nabla_p\nabla^k f)$$

and
$$\tilde{\Gamma}^p_{00} = -(\tilde{R}^p_0 - \tilde{\nabla}_0\tilde{\nabla}^p f) = \nabla^p(-\frac{1}{2}R + \frac{\partial f}{\partial t} - \frac{1}{2}|\nabla f|^2)$$

to obtain the last equality. This completes the proof of the theorem.

Similar to the previous section, if $(\tilde{g}_\tau, \tilde{\nabla}_\tau)$ satisfies the modified Ricci flow for degenerate metrics, then the evolution of the Riemann curvature operator is given by:
$$\frac{\partial}{\partial t}\tilde{R}m = \tilde{\Delta}\tilde{R}m + \mathcal{L}_{\tilde{\nabla} f}\tilde{R}m + \tilde{R}m^2 + \tilde{R}m^\#.$$

# 4  An Approximation approach

In this section we consider a two-parameter family of Riemannian (nondegenerate) metrics $\tilde{g}_{\epsilon,\delta}$ on $\tilde{M}_0 = M \times [0,t)$ and obtain the connection $\tilde{\nabla}_\tau$ defined by (A1)-(A4) as the limit of the Levi-Civita connections of $\tilde{g}_{\epsilon,\delta}$ as $\epsilon$ and $\delta$ tend to infinity[7]. On the other hand, the Harnack quantity $Z$ is the limit of the Riemann curvature tensors of $\tilde{g}_{\epsilon,\delta}$ as $\epsilon$ tends to infinity and $\delta$ tends to zero. We define the metrics by:

$$\tilde{g}_\epsilon(x,t) = g(x,t) + \left(R + \frac{\epsilon}{2(t+\delta)}\right) dt^2.$$

This metric is positive-definite at points where $R + \epsilon(2t+2\delta)^{-1} > 0$. In particular, for any $\delta > 0$ $\tilde{g}_{\epsilon,\delta}$ is positive-definite on compact subsets of $\tilde{M}_0$ provided $\epsilon$ is large enough.

**Remark 4.1** *The metric $\tilde{g}_{\epsilon,\delta}$ induces a metric on $\wedge^2 T^* \tilde{M}$. A two-form on $\tilde{M}_0$ is a section in $\wedge^2 T^* \tilde{M}$, and we have*

$$\wedge^2 T^*_{(x,t)} \tilde{M} \equiv \wedge^2(T_x^* M \oplus R_t) \cong \wedge^2 T_x^* M \oplus \wedge^1 T_x^* M \text{ as vector bundles.}$$

*The second identification is given by*

$$\begin{aligned} dx^\alpha \wedge dx^\beta &\to dx^\alpha \wedge dx^\beta \\ dx^\gamma \wedge dt &\to \frac{1}{\sqrt{R + \frac{\epsilon}{2(t+\delta)}}} dx^\gamma. \end{aligned}$$

*As $\epsilon$ tends to infinity, we get the semi-direct Lie algebra structure on $\wedge^2 T^* \tilde{M}$ given by (1.5) and the degenerate metric on $\wedge^2 T^* \tilde{M}$ given by (1.6).*

**Remark 4.2** $\tilde{g}_{ij} = g_{ij}$, $\tilde{g}^{ij} = g^{ij}$ and $\tilde{g}_{0i} = \tilde{g}^{0i} = 0$ for $1 \leq i,j \leq n$. $\tilde{g}_{00} = R + \epsilon(2t+2\delta)^{-1}$.

We now study the asymptotic behavior of the connection and curvature operator of $(\tilde{M}^{n+1}, g_{\epsilon,\delta})$ as $\epsilon \to \infty$.

**Lemma 4.3** *For all $\epsilon, \delta > 0$, the Levi-Civita connections of $\tilde{g}_{\epsilon,\delta}$ are given by:*

1) $\tilde{\Gamma}^k_{ij} = \Gamma^k_{ij}$ for $1 \leq i,j,k \leq n$,
2) $\tilde{\Gamma}^0_{ij} = \frac{R_{ij}}{R + \epsilon(2t+2\delta)^{-1}}$ for $1 \leq i,j \leq n$,
3) $\tilde{\Gamma}^k_{i0} = -R^k_i$ for $1 \leq i,k \leq n$,
4) $\tilde{\Gamma}^k_{00} = -\frac{1}{2} \nabla^k R$ for $1 \leq k \leq n$,
5) $\tilde{\Gamma}^0_{i0} = \frac{\nabla_i R}{2(R + \epsilon(2t+2\delta)^{-1})}$ for $1 \leq i \leq n$,
6) $\tilde{\Gamma}^0_{00} = \frac{\partial_t R}{2(R + \epsilon(2t+2\delta)^{-1})} + \frac{-\epsilon}{4t^2(R + \epsilon(2t+2\delta)^{-1})}$.

*at points where $\tilde{g}_{\epsilon,\delta}$ is positive-definite.*

---

[7] If we just let $\epsilon$ tend to infinity, the connection $\tilde{\nabla}_\tau$ is the limit except for the component $\tilde{\Gamma}^0_{00}$ (see formula 6 of Lemma 4.1) which tends to $-(2t+2\delta)^{-1}$ as $\epsilon \to \infty$, whereas by (A2), $\tilde{\Gamma}^0_{00} = 0$. Of course, this discrepancy disappears if we let $\delta \to \infty$.

Proof. Recall the Christoffel symbols are given by

$$\tilde{\Gamma}^k_{ij} = \frac{1}{2}\tilde{g}^{kl}(\partial_i \tilde{g}_{jl} + \partial_j \tilde{g}_{il} - \partial_l \tilde{g}_{ij}) \text{ for } 0 \leq i, j, k, l \leq n.$$

The lemma follows from this and a straightforward computation using Remark 4.2.

The Riemann curvature tensors of the metrics $\tilde{g}_{\epsilon,\delta}$ are given by

**Lemma 4.4**

1) $\tilde{R}^l_{ijk} = R^l_{ijk} - (R^l_i R_{jk} - R^l_j R_{ik})(R + \epsilon(2t + 2\delta)^{-1})^{-1}$ for $1 \leq i, j, k, l \leq n$,
2) $\tilde{R}^l_{0jk} = -\nabla_k R^l_j + \nabla^l R_{jk} - \frac{1}{2}(R_{jk}\nabla^l R - R^l_j \nabla_k R)(R + \epsilon(2t + 2\delta)^{-1})^{-1}$ for $1 \leq j, k, l \leq n$,
3) $\tilde{R}^l_{i00} = \partial_t R^l_i - \frac{1}{2}\nabla_i \nabla^l R - R^m_i R^l_m - [(\frac{\partial_t R}{2} + \frac{-\epsilon}{4t^2})R^l_i - \frac{\nabla_i R \nabla^l R}{4}](R + \epsilon(2t + 2\delta)^{-1})^{-1}$ for $1 \leq i, l \leq n$.

Proof.

1. For $1 \leq i, j, k, l \leq n$, we have

$$\begin{aligned}\tilde{R}^l_{ijk} &= \partial_i \tilde{\Gamma}^l_{jk} - \partial_j \tilde{\Gamma}^l_{ik} + \tilde{\Gamma}^m_{jk}\tilde{\Gamma}^l_{im} - \tilde{\Gamma}^m_{ik}\tilde{\Gamma}^l_{jm} \\ &= \partial_i \Gamma^l_{jk} - \partial_j \Gamma^l_{ik} + \sum_{1 \leq m \leq n}(\Gamma^m_{jk}\Gamma^l_{im} - \Gamma^m_{ik}\Gamma^l_{jm}) + \tilde{\Gamma}^0_{jk}\tilde{\Gamma}^l_{i0} - \tilde{\Gamma}^0_{ik}\tilde{\Gamma}^l_{j0} \\ &= R^l_{ijk} - (R^l_i R_{jk} - R^l_j R_{ik})(R + \epsilon(2t + 2\delta)^{-1})^{-1}\end{aligned}$$

2. For $1 \leq l, j, k \leq n$, we have

$$\begin{aligned}\tilde{R}^l_{0jk} &= \partial_0 \tilde{\Gamma}^l_{jk} - \partial_j \tilde{\Gamma}^l_{0k} + \tilde{\Gamma}^m_{jk}\tilde{\Gamma}^l_{0m} - \tilde{\Gamma}^m_{0k}\tilde{\Gamma}^l_{jm} \\ &= \partial_0 \Gamma^l_{jk} - \partial_j \Gamma^l_{k0} + \sum_{1 \leq m \leq n}(-\Gamma^m_{jk} R^l_m + \Gamma^l_{jm} R^m_k) + \tilde{\Gamma}^0_{jk}\tilde{\Gamma}^l_{00} - \tilde{\Gamma}^0_{0k}\tilde{\Gamma}^l_{j0} \\ &= -(\nabla_j R^l_k + \nabla_k R^l_j - \nabla^l R_{jk}) + \nabla_j R^l_k - \frac{1}{2}(R_{jk}\nabla^l R - R^l_j \nabla_k R)(R + \epsilon(2t + 2\delta)^{-1})^{-1}\end{aligned}$$

3. For $1 \leq i, l \leq n$, we have

$$\begin{aligned}\tilde{R}^l_{i00} &= \partial_i \tilde{\Gamma}^l_{00} - \partial_0 \tilde{\Gamma}^l_{i0} + \tilde{\Gamma}^m_{00}\tilde{\Gamma}^l_{im} - \tilde{\Gamma}^m_{i0}\tilde{\Gamma}^l_{0m} \\ &= \partial_i(-\frac{1}{2}\nabla^l R) - \partial_0(-R^l_i) + \sum_{1 \leq m \leq n}(-\frac{1}{2}\Gamma^l_{im}\nabla^m R - R^m_i R^l_m) + \tilde{\Gamma}^0_{00}\tilde{\Gamma}^l_{i0} - \tilde{\Gamma}^0_{i0}\tilde{\Gamma}^l_{00} \\ &= \partial_0 R^l_i - \frac{1}{2}\nabla_i\nabla^l R - R^m_i R^l_m + \tilde{\Gamma}^0_{00}\tilde{\Gamma}^l_{i0} - \tilde{\Gamma}^0_{i0}\tilde{\Gamma}^l_{00} \\ &= \partial_t R^l_i - \frac{1}{2}\nabla_i\nabla^l R - R^m_i R^l_m - [(\frac{\partial_t R}{2} + \frac{-\epsilon}{4t^2})R^l_i - \frac{\nabla_i R \nabla^l R}{4}](R + \epsilon(2t + 2\delta)^{-1})^{-1}\end{aligned}$$

Taking the appropriate limits, we obtain the connection $\tilde{\nabla}_\tau$ and the Harnack quantity.

**Theorem 4.5** 1. The connection $\tilde{\nabla}_\tau$ is the limit of $\tilde{\nabla}_{\epsilon,\delta}$ as $\epsilon, \delta \to \infty$.

2. Hamilton's Harnack quantity $Z + \frac{1}{2t}R_{ij}W_i W_j$ is the limit of

$$\tilde{Rm}_{\epsilon,\delta}(U \oplus W, U \oplus W) \text{ as } \epsilon \to \infty \text{ and } \delta \to 0.$$

**Acknowledgement.** We would like to thank Richard Hamilton for very helpful suggestions and discussions.